\documentclass[11pt]{article}
\usepackage{amssymb,amsbsy,amsmath,amsfonts,amssymb,amscd}
\usepackage{latexsym}

\newtheorem{lemma}{Lemma}[section]
\newtheorem{theorem}[lemma]{Theorem}
\newtheorem{corollary}[lemma]{Corollary}
\newtheorem{proposition}[lemma]{Proposition}
\newtheorem{definition}[lemma]{Definition}

\begin{document}

\def\C{{\mathbb C}}
\def\N{{\mathbb N}}
\def\Z{{\mathbb Z}}
\def\R{{\mathbb R}}
\def\PP{\cal P}
\def\p{\rho}
\def\phi{\varphi}
\def\ee{\epsilon}
\def\ll{\lambda}
\def\l{\lambda}
\def\a{\alpha}
\def\b{\beta}
\def\D{\Delta}
\def\g{\gamma}
\def\rk{\text{\rm rk}\,}
\def\dim{\text{\rm dim}\,}
\def\ker{\text{\rm ker}\,}
\def\square{\vrule height6pt width6pt depth 0pt}
\def\epsilon{\varepsilon}
\def\phi{\varphi}
\def\kappa{\varkappa}
\def\wz{\thinspace}
\def\proof{P\wz r\wz o\wz o\wz f.\hskip 6pt}
\def\quest#1{\hskip5pt {\scshape  Problem} {\rm #1}.\hskip 6pt}
\def\leq{\leqslant}
\def\geq{\geqslant}
\def\pd#1#2{\frac{\partial#1}{\partial#2}}
\def\limsup{\mathop{\overline{\hbox{\rm lim}\,}}}
\def\ug#1#2{\left\langle#1,#2\right\rangle}
\def\kk#1#2{{\k\langle#1,#2\rangle}}
\def\sv{\bf{ k} \langle X \rangle}
\def\k{k }
\def\lxr{\langle X \rangle}
\def\defin#1{\smallskip\noindent
{\scshape  Definition} {\bf #1}{\bf .}\hskip 8pt\sl}

\def\doubarr#1#2{\mathop{\hbox{$\vcenter{\offinterlineskip\halign
{\kern2pt\hfil##\hfil\kern2pt\cr \vrule height6pt width0pt
depth0pt\cr \smash{${\longleftarrow}\!\!{-}\!\!{-}$}\cr \vrule
height4pt width0pt depth0pt\cr
\smash{${\longleftarrow}\!\!{-}\!\!{-}$}\cr
}}$}}\limits_{#1}^{#2}}

\def\bull#1{\mathop{\bullet}\limits_{#1}}

\def\siglearr#1{\mathop{\hbox{${-}\!\!{-}\!\!{\longrightarrow}$}}\limits^{#1}}

\font\LINE=line10 scaled 1440 \font\Line=line10 \font\bIg=cmr10
scaled 1728 \font\biG=cmr10 scaled 1440
\def\lini{{\LINE\char"40}}
\def\li{{\Line\char"40}}
\def\pph{\vrule height4pt width0pt depth0pt}
\def\pha{{\phantom0}}
\def\sst{\scriptstyle}
\def\ddd{\displaystyle}
\def\>{\multispan}

\def\ramka#1#2#3{\hbox{$\vcenter{\offinterlineskip\halign
{\vrule\vrule\kern#2pt\hfil##\hfil\kern#2pt\vrule\vrule\cr
\noalign{\hrule}
\noalign{\hrule}
\vrule height #3pt depth0pt width0pt\cr
#1\cr
\vrule height #3pt depth0pt width0pt\cr
\noalign{\hrule}
\noalign{\hrule}}
}$}}

\def\jdots{\hbox{$\vcenter{\offinterlineskip\halign
{\hfil##\hfil\kern2pt&\hfil##\hfil\kern2pt&\hfil##\hfil\cr
&&.\cr
\noalign{\kern3pt}
&.&\cr
\noalign{\kern3pt}
.&&\cr
}}$}}

\def\Mone{\,\,\,\,\,\,\hbox{$\vcenter{\offinterlineskip\halign
{##\hfil&##\hfil&##\hfil&##\hfil&##\hfil&##\hfil&##\hfil\cr
&&&&&\smash{$\sst\bullet$}&$\sst n-1$\cr
\lini&\lini&\lini&&&\smash{$\pha\atop\ddd\vdots$}&\cr
\lini&\lini&\lini&\lini&&&\cr
\lini&\lini&\lini&\lini&\lini&&\smash{${}_2$}\cr
&\lini&\lini&\lini&\lini&&\smash{${}_1$}\cr
&&\lini&\lini&\lini&&\smash{${}_0$}\cr
\smash{$\sst\bullet$}&\rlap{\ $\dots$}&&&&&\cr
\vrule height 3pt depth0pt width0pt\cr&&&&&&\cr
\llap{$\sst-n$}$\sst+1$&&\rlap{$\,\,\,\sst-2$}&\rlap{$\,\,\,\sst-1$}&&&\cr
}}$}\,}

\def\MtwoA{\hbox{$\vcenter{\offinterlineskip\halign
{##\hfil&##\hfil&##\hfil&##\hfil&##\hfil\cr
0&\ \ 1&&\ \smash{\lower 6pt\hbox{\biG0}}&\cr
\noalign{\kern-2pt}
&$\,\ddots$&$\ddots$&&\cr
\noalign{\kern-6pt}
&&$\ddots$&$\ddots$&\cr
\noalign{\kern-6pt}
&\smash{\biG0}&&$\ddots$&\smash{\raise 6pt\hbox{\ 1}}\cr
\noalign{\kern6pt}
&&&&\smash{\raise 1pt\hbox{\ 0}}\cr
}}$}}

\def\MtwoB{\,\hbox{$\vcenter{\offinterlineskip\halign
{##\hfil&##\hfil&##\hfil&##\hfil&##\hfil&##\hfil\cr
&&&\rlap{\large$\dots$}&&\smash{$\sst\bullet$}\cr
\lini&\rlap{\lini}\kern.5pt\rlap{\lini}\kern.5pt\lini\kern-1pt&\lini&&&{\large\smash{$\pha\atop\ddd\vdots$}}\cr
&\lini&\rlap{\lini}\kern.5pt\rlap{\lini}\kern.5pt\lini\kern-1pt&\lini&&\cr
&&\lini&\rlap{\lini}\kern.5pt\rlap{\lini}\kern.5pt\lini\kern-1pt&\lini&\cr
&\smash{\bIg0}&&\lini&\rlap{\lini}\kern.5pt\rlap{\lini}\kern.5pt\lini\kern-1pt&\cr
&&&&\lini&\cr
}}$}}

\def\Mfour{\,\hbox{$\vcenter{\offinterlineskip\halign
{\hfil##&##&##&##\hfil\cr
\ramka{$J_1$}{7}{7}&&\ \bIg0\cr
&\ramka{$J_2$}{10}{10}&&\cr
\noalign{\kern-5pt}
\smash{\lower 5pt\hbox{\bIg0\ }}&&$\ddots$&\cr
&&&\ramka{$J_m$}{4}{5}\cr
}}$}\,}

\def\Mfive{\,\hbox{$\vcenter{\offinterlineskip\halign
{\vrule\hfil##\hfil\vrule&\hfil##\hfil\vrule&\hfil##\hfil\vrule&\hfil##\hfil\vrule\cr
\noalign{\hrule}
\ramka{$A_{11}$}{9}{12}&$A_{12}$&&$A_{1m}$\cr
\noalign{\hrule}
&\ramka{$A_{22}$}{12}{15}&&\cr
\noalign{\hrule}
\noalign{\kern-5pt}
&&$\ddots$&\cr
\noalign{\hrule}
$A_{m1}$&$A_{m2}$&&\ramka{$A_{mm}$}{4}{10}\cr
\noalign{\hrule}
}}$}\,}

\def\Msix{\,\hbox{$\vcenter{\offinterlineskip\halign
{\vrule\hfil##\hfil\vrule&\hfil##\hfil\vrule&\hfil##\hfil\vrule&\hfil##\hfil\vrule\cr
\noalign{\hrule}
\ramka{$[A_{11}]$}{9}{14}&$[A_{12}]$&&$[A_{1m}]$\cr
\noalign{\hrule}
&\ramka{$[A_{22}]$}{12}{17}&&\cr
\noalign{\hrule}
\noalign{\kern-5pt}
&&$\ddots$&\cr
\noalign{\hrule}
$[A_{m1}]$&$[A_{m2}]$&&\ramka{$[A_{mm}]$}{4}{12}\cr
\noalign{\hrule}
}}$}\,}

\vfill\break

\def\MsevenA{\,\hbox{$\vcenter{\offinterlineskip\halign
{\vrule\kern65pt##\hfil&##\hfil&##\hfil&##\hfil&##\hfil&##\hfil\vrule\cr
&&&&&\rlap{\smash{\lower2pt\hbox{\kern-2pt$\sst\bullet$}}}\cr
\noalign{\hrule}
\lini&\lini&\lini&&\rlap{\kern-4pt\smash{$\jdots$}}&\cr
&\lini&\lini&\lini&&\cr
\kern-30pt\smash{\bIg0}&&\lini&\lini&\lini&\cr
&&&\lini&\lini&\cr
&&&&\lini&\cr
\noalign{\hrule}
}}$}\,}

\def\MsevenB{\,\hbox{$\vcenter{\offinterlineskip\halign
{\vrule##\hfil&##\hfil&##\hfil&##\hfil&##\hfil&##\hfil\vrule\cr
&&&&&\rlap{\smash{\lower2pt\hbox{\kern-2pt$\sst\bullet$}}}\cr
\noalign{\hrule}
\lini&\lini&\lini&&\rlap{\kern-4pt\smash{$\jdots$}}&\cr
&\lini&\lini&\lini&&\cr
&&\lini&\lini&\lini&\cr
&&&\lini&\lini&\cr
&&\smash{\lower30pt\hbox{\bIg0}}&&\lini&\cr
\vrule height65pt depth0pt width0pt &&&&&\cr
\noalign{\hrule}
}}$}\,}

\def\Meighta{\hbox{$\vcenter{\offinterlineskip\halign
{\vrule\vrule##\hfil&##\hfil&##\hfil&##\hfil\vrule\vrule\cr
\noalign{\hrule}
\noalign{\hrule}
\li&\rlap{\li}\kern.4pt\rlap{\li}\kern.4pt\li\kern-.8pt&\li&\cr
&\li&\rlap{\li}\kern.4pt\rlap{\li}\kern.4pt\li\kern-.8pt&\li\cr
&&\li&\rlap{\li}\kern.4pt\rlap{\li}\kern.4pt\li\kern-.8pt\cr
&&&\li\cr
\noalign{\hrule}
\noalign{\hrule}
}}$}}

\def\Meightaa{\hbox{$\vcenter{\offinterlineskip\halign
{##\hfil&##\hfil&##\hfil&##\hfil\cr
\li&\li&\li&\cr
&\li&\li&\li\cr
&&\li&\li\cr
&&&\li\cr
}}$}}

\def\Meightb{\hbox{$\vcenter{\offinterlineskip\halign
{\vrule\vrule##\hfil&##\hfil&##\hfil&##\hfil&##\hfil\vrule\vrule\cr
\noalign{\hrule}
\noalign{\hrule}
\li&\rlap{\li}\kern.4pt\rlap{\li}\kern.4pt\li\kern-.8pt&\li&&\cr
&\li&\rlap{\li}\kern.4pt\rlap{\li}\kern.4pt\li\kern-.8pt&\li&\cr
&&\li&\rlap{\li}\kern.4pt\rlap{\li}\kern.4pt\li\kern-.8pt&\li\cr
&&&\li&\rlap{\li}\kern.4pt\rlap{\li}\kern.4pt\li\kern-.8pt\cr
&&&&\li\cr
\noalign{\hrule}
\noalign{\hrule}
}}$}}

\def\Meightc{\hbox{$\vcenter{\offinterlineskip\halign
{\vrule\vrule##\hfil&##\hfil&##\hfil\vrule\vrule\cr
\noalign{\hrule}
\noalign{\hrule}
\li&\rlap{\li}\kern.4pt\rlap{\li}\kern.4pt\li\kern-.8pt&\li\cr
&\li&\rlap{\li}\kern.4pt\rlap{\li}\kern.4pt\li\kern-.8pt\cr
&&\li\cr
\noalign{\hrule}
\noalign{\hrule}
}}$}}

\def\Meightcc{\hbox{$\vcenter{\offinterlineskip\halign
{##\hfil&##\hfil&##\hfil\cr
\li&\li&\li\cr
&\li&\li\cr
&&\li\cr
}}$}}

\def\Meight{\,\,\hbox{$\vcenter{\offinterlineskip\halign
{\vrule\hfil##\vrule&\hfil##\vrule&\hfil##\hfil\vrule&\hfil##\vrule\cr
\noalign{\hrule}
\Meighta\llap{\lower7pt\hbox{0}\kern23pt}&\Meightaa\llap{\lower7pt\hbox{0}%
\kern33pt}&&\raise5.5pt\Meightcc\llap{\lower7pt\hbox{0}\kern18pt}\cr
\noalign{\hrule}
\raise5pt\Meightaa\llap{\lower7pt\hbox{0}\kern23pt}&\Meightb\llap{\lower7pt%
\hbox{0}\kern33pt}&&\raise10pt\Meightcc\llap{\lower7pt\hbox{0}\kern18pt}\cr
\noalign{\hrule}
&&\raise5pt\hbox{$\,\ddots\,$}&\cr
\noalign{\hrule}
\Meightcc\llap{\lower5pt\hbox{0}\kern23pt}&\Meightcc\llap{\lower5pt%
\hbox{0}\kern33pt}&&\Meightc\llap{\lower5pt\hbox{0}\kern18pt}\cr
\noalign{\hrule}
}}$}\,\,}


\vskip1cm

\title{The representation theory of the Jordanian algebra}

\author{
 Natalia K. Iyudu}

\date{}

\maketitle

\small

\centerline{Queen Mary, University of London, Mile End Road, London
E1 4NS U.K., }

\smallskip

\centerline{Department of Pure Mathematics, Queen's University
Belfast, Belfast BT7 1NN, U.K.}

\smallskip

\centerline{ {\bf e-mail:}  \,\, n.iyudu@qub.ac.uk}

\bigskip

\small

{\bf Abstract}

\bigskip

We describe the complete set of pairwise non-isomorphic irreducible
modules $S_{\a}$ over the algebra $R=k\langle x,y\rangle/
(xy-yx-y^2)$, and the rule how they could be glued to
indecomposables. Namely, we show that  ${\rm
Ext}_k^1(S_{\a},S_{\b})=0$, if $\a \neq \b $. Also the set of all
representations  is described subject to the Jordan normal form of
$Y$.

We study then properties of the image algebras in the endomorphism
ring. Among facts we prove is that they are all basic algebras.
Along this line we establish an analogue of the
Gerstenhaber--Taussky--Motzkin theorem on the dimension of algebras
generated by two commuting matrices. All image algebras of
indecomposable modules turned out to be local complete algebras. We
compare them with the Ringel's classification by means of finding
relations of image algebras. As a result we derive that all image
algebras of $n$-dimensional representations with full block $Y$ are
tame for $n \leq 4$ and wild for $m \geq 5$.

We suggest a stratification of representation space of $R$ by
 partitions of $n$ related to the Jordan normal form
of $Y$. We give a complete classification by parameters for some
strata and present examples of tame and even finite type (up to
automorphisms) strata, while the generic stratum is wild.


\normalsize

\bigskip
\tableofcontents
\medskip

\section{Introduction}

We consider here a quadratic algebra given by the following
presentation: $R=k\langle x,y\rangle/ (xy-yx-y^2)$. This algebra
appeared in various different contexts in mathematics and physics.
First of all it is a kind of a quantum plane: one of the two
Auslander regular algebras of global dimension two in the
Artin--Shelter classification \cite{AS}. The other one is a usual
quantum plane $k\langle x,y\rangle/ (xy-qyx)$. There were studied,
for example, deformations of $GL(2)$ analogues to $GL_q(2)$ with
respect to this algebra in 80-90th in Manin's 'Quantum group'
\cite{Ma}, \cite{Ko},
 where this algebra appeared under the name
Jordanian algebra.

This algebra is  also a simplest element in the class of  RIT
(relativistic internal time) algebras. The latter appeared  and
investigations were started in papers \cite{a1}, \cite{ia},
\cite{wia}, \cite{iw}, \cite{cs}. The class of RIT algebras arises
from a modification of the Poincare algebra of the Lorenz group
SO(3,1) by means of introducing the additional generator
corresponding to the relativistic internal time. The algebra $R$
above is a RIT algebra of type (1,1). Our studies of this algebra
reported in \cite{MP}, \cite{td}.

 Let we mention
that  algebra $R$ is a subalgebra of the first Weyl algebra $A_1$.
The latter has no finite dimensional representations, but $R$ turned
out to have quite a rich structure of them. Category of finite
dimensional modules over $R$ contains, for example, as a full
subcategory ${\rm mod} GP(n,2)$, where $GP(n,2)$ is a
Gelfand--Ponomarev algebra \cite{GP} with the nilpotency degrees of
variables $x$ and $y$,  $n$ and $2$ respectively. On the other hand
we show in section \ref{rfd} that $R$ is residually finite
dimensional.

 We are interested here in
representations over algebraically closed field $k$ of
characteristic $0$. In few places we suppose $k= \Bbb C$, this will
be pointed out separately. We denote throughout be the category of
all $R$-modules by $\rm{Mod }R$, the category of finite dimensional
$R$-modules by ${\rm mod} R$ and $\rho_n \in {\rm mod} R$ stands for
an $n$-dimensional representation of $R$.

We study here the category ${\rm mod} R$ of finite dimensional
representations of infinite dimensional algebra $R$ first by
encountering some properties of
 its finite dimensional images in the endomorphism ring
${\rm End} (k^n)$. Then we  suggest a stratification of a
representation space related to partitions which define the Jordan
structure of $Y=\rho_n(y)$ and give a classification and tameness
results for some strata.


Toward the first approach we prove (section 2) that images in
endomorphism rings are {\it basic} algebras, that is their
semisimple parts are a direct sums of fields. This allow to
associate a {\it quiver} to any representation and to classify
representations using these quivers. It turns out that
indecomposable modules have usually a typical wild quiver with one
vertex and two loops and in some cases the quiver with one vertex
and one loop. The simple, but important fact for the structural
properties of image algebras is that  $Y=\rho_n(y)$ is nilpotent for
any $\rho_n \in mod R$. Note that this is not necessarily   the case
when the characteristic of the basic field is not zero. After the
description of all finite dimensional modules subject to Jordan
normal form of $Y$ in section \ref{descript}, we study {\it
irreducible} and {\it indecomposable} modules in sections
\ref{ind},\ref{eqc}. We describe the complete set of irreducible
modules $S_{\alpha}$ and show how one could glue them together:
${\rm Ext}_k^1(S_{\a},S_{\b})=0$, if $\a \neq \b $. This means that
indecomposables have always $X$ with only one eigenvalue.

From the above results we see that any algebra which is an image of
indecomposable representation $\p_n$ is a {\it local complete}
algebra. Hence we could apply Ringel's classification \cite{Ringloc}
of local complete algebras to those images and after calculating
defining relations of image algebras get that all of them are {\it
tame} for $n \leq 4$ and for $n \geq 5$ they are {\it wild}. These
results are described in section \ref{ring}.

In section \ref{Ger} we prove an analogue of the
Gerstenhaber--Taussky--Motzkin theorem \cite{Ger}, \cite{Gur} on the
dimension of algebras generated by two commuting matrices. The
dimension of image algebras of representations of $R$ does not
exceed $n(n+2)/4$ for even $n$ and $(n+1)^{2} /4$ for odd $n$. This
estimate is attained  for the family of representations with the
full block Jordan form of $Y$.

In section \ref{param} we consider an action of $GL_n$ on the
representation space of the algebra $R$: {\it mod(R,n)}, which is a
representation space of a {\it wild quiver with relation}. We take
as a model strata the one related to the full Jordan block $Y$ and
show that it is parametrizable by two parameters in a conventional
sense, and has a finite representation type with respect to
auto-equivalence relation on reps (defined via gluing of isoclasses,
which are coincide up to automorphisms of initial algebra, see
section \ref{rfd} for a precise definition). To get results on the
auto-equivalence we describe the automorphism group of $R$ in
section \ref{aut}. We give  examples of {\it tame strata} (in proper
sense, e.i.   not of finite representation type) up to
auto-equivalence. Let we mention, that normally, for generic
partition, the stratum is wild for this algebra.

Main tool we use for the parametrization results is to consider in
stead of the whole action of $GL_n$ on the representation space the
action of the centralizer of Jordan form of $Y$ on those points of
the space where $Y$ is in this Jordan form. While the group which
acts is not reductive any more, the space where it acts become much
simpler. Due to 1-1 correspondence of orbits under these two actions
one can lift results on classification  from one setting to another.

\section{Structural properties of the images of
representations and quivers}\label{quiv}

We consider here the case when $k$ is an algebraically closed field
of characteristic zero. Let $\p: R \rightarrow {\rm End} (k^n)$ be
an arbitrary finite dimensional representation of $R$, denote by
$A_{\p,n}=\p (R)$ an image of $R$ in the endomorphism ring. We will
write also $A_n$ or $A$ when it is clear from the context which $\p$
and $n$ we mean.

We derive in this section some structural properties   of algebras
$A_{n,\p}$. They all turned out to be basic; in any of such algebra
the image of $y$ is nilpotent; complete system of orthogonal
idempotents in $A$ corresponds to the different eigenvalues of $x$,
etc.

Let $J(A)=J$ be the Jacobson radical of the algebra $A_{\p,n}$. We
show first that $A$ is {\it basic}, that is its semisimple part
$A/J(A)$ is a direct product of division rings, or in our case of
algebraically closed field $k$ it is the same as direct product of
several copies of the field $k$. Basic algebras take their name
particularly because they are basic from the point of view of Morita
equivalence. Due to the Wedderburn--Artin theorem and the
equivalence of categories of modules Mod-$R$ and Mod-$R_n$ ($R_n$
are $n\times n$ matrices over $R$), any artinian semisimple algebra
is Morita equivalent to a finite direct sum of division rings. So we
show that image algebras of all finite dimensional representations
of $R$ has this special place between all finite dimensional
algebras in sense of Morita equivalence. We also prove in section
\ref{rfd} that $R$ is residually finite dimensional, and together
with the above fact it gives as a consequence that $R$ is residually
basic.

Lemmata below describes the  structural properties of image algebras
for $R$. The following fact probably allow many different proofs. We
present here the shortest we know.

\begin{lemma}\label{l2.1} Let $Y=\rho_n(y)$. Then the matrix $Y$ is nilpotent.
\end{lemma}

\proof Suppose that the matrix $Y$ is not nilpotent and hence has a
nonzero eigenvalue. We take a projector $P$ on the subspace
corresponding to this eigenvalue. It is obviously commute with any
matrix, particularly with $Y: PY=YP$, and is an idempotent operator:
$P^2=P$. Hence multiplying our relation $XY-YX=Y^2$ from the right
and from the left hand  side by $P$ and using above two notices we
can observe that operators $X'=PXP$ and $Y'=PYP$ also satisfy the
same relation: $X'Y'-Y'X'={Y'}^2$. Taking into account that $Y'$ has
a form of one or more Jordan blocks with the same nonzero eigenvalue
$\ll$, we get that traces of right and left parts of the relation
can not coincide. This contradiction complete the proof. $\Box$

Let we prove here also a little bit more general fact.

\begin{lemma} \label{nilnil}
Let $X,Y$ be $n\times n$ matrices over an algebraically closed field
$k$ of characteristic zero. Assume that the commutator $Z=XY-YX$
commutes with $Y$. Then $Z$ is nilpotent.
\end{lemma}

\proof
Assume the contrary.
Then $Z$ has a non-zero eigenvalue in $z\in k$. Let
$$
L={\rm ker}\,(Z-zI)^n\ \ \text{and}\ \ N={\rm Im}(Z-zI)^n.
$$
The subspace $L$ is known as a main subspace for $Z$ corresponding
to the eigenvalue $z$. Clearly $L\neq \{0\}$. It is well-known that
$k^n$ is the direct sum $k^n=L\oplus N$ of $Z$-invariant linear
subspaces $L$ and $N$. Due to $ZY=YZ$, the subspaces $L$ and $N$ are
also invariant for $Y$.

Consider the linear projection $P$ along $N$ onto $L$. Since $L$ and
$N$ are invariant under both $Y$ and $Z$, we have $ZP=PZ$ and
$YP=PY$. Multiplying the equality  $Z=XY-YX$ by $P$ from the left
and from the right hand side and using the equalities $ZP=PZ$,
$YP=PY$ and $P^2=P$, we get
$$
ZP=PXPY-YPXP.
$$
Since $ZP$ vanishes on $N$ and $ZP-zI$ has only one eigenvalue $z$,
then after restriction to $L$, we have ${\rm tr}\,ZP=z\,{\rm
dim}\,L$. On the other hand ${\rm tr}\,PXPY={\rm tr}\,YPXP$ since
the trace of a product of two matrices does not depend on the order
of the product. Thus, the last display implies that
$
z\,{\rm dim}\,L=0,
$
which is not possible since $z\neq 0$ and ${\rm dim}\,L>0$. $\Box$

Coming back to the case of simplest RIT algebra, we have further

\begin{lemma}\label{l2.2} Let $X=\rho_n(x)$. Then the matrix
$S=S(X)=(X-\lambda_1 I) \dots (X-\lambda_r I)$ is nilpotent.
\end{lemma}

\proof Note that Spec$\,p(X)=p({\rm Spec}\,X)$ for any polynomial
$p$. Spec$\,X$ in our case is $\{\lambda_1,\dots,\lambda_r\}$ and
hence Spec$\,S=\{0\}$. Therefore the matrix $S$ is nilpotent. $\Box$

\begin{lemma}\label{l2.3} Any nilpotent element of the algebra $A=\rho(R)$
belongs to the radical $J(A)$.  \end{lemma}

\proof We will use here the feature of an algebra $A$ that it has
the presentation as a quotient of free algebra containing our main
relation. Namely, it has a presentation: $A=\k\langle x,y|
xy-yx=y^2, R_A \rangle$, where $R_A \subset \k\langle x,y \rangle$
is the set of additional relations specific for the given image
algebra. Thus we can think of elements in $A$ as of polynomials in
two variables (subject to some relations). Let $Q(x)$ be a
polynomial on one variable $Q(x)\in \k[x]$ and $Q(X)\in A$ be a
nilpotent element with the degree of nilpotency $N$: $Q^N=0$. We
show first that $Q\in J(A)$. We have to check that for any
polynomial $a\in \k\langle x,y\rangle$, $1-a(X,Y)Q(X)$ is
invertible. It suffices to verify that $a(X,Y)Q(X)$ is nilpotent. By
lemma \ref{l2.1} $Y$ is nilpotent. Denote by $m$ the degree of
nilpotency of $Y$: $Y^m=0$. Let we verify that
$(a(X,Y)Q(X))^{mN}=0$. Present $a(X,Y)$ as $u(X)+Yb(X,Y)$. If then
we consider a word of length not less then $mN$ on letters
$\alpha=u(X)Q(X)$ and $\beta=Yb(X,Y)Q(X)$ then we can see that it is
equal to zero. Indeed, if there are at least $m$ letters $\beta$
then using the relation $XY-YX=Y^2$ to commute the variables one can
rewrite the word as a sum of words having a subword $Y^m$. Otherwise
our word has the subword $\alpha^N=u(X)^NQ(X)^N=0$. Thus, $Q(X)\in
J(A)$.

Note now that if we have an arbitrary nilpotent polynomial
 $G(X,Y)$, we can separate the terms containing $Y$:
 $G(X,Y)=Q(X)+YH(X,Y)$. To obtain nilpotency of any element
 $a(X,Y)G(X,Y)$ it suffices to verify nilpotency of
 $a(X,Y)Q(X)$, which was already proven, because the
 relation $[X,Y]=Y^2$ allows to commute with $Y$,
 preserving (or increasing) the degree of $Y$. $\Box$

\begin{corollary}\label{c2.1} The Jacobson radical of $A=\rho(R)$ consists
precisely of all nilpotent elements. \end{corollary}

Particularly,

\kern-3mm

\begin{corollary}\label{c2.2} Let $Y=\rho(y)$. Then $Y\in J(A)$. \end{corollary}

Let we formulate here also another property of the radical, which
will be on use later on.

\begin{corollary}\label{radft}
The Jacobson radical of $A=\rho(R)$ consists of all polynomials on
$X=\p(x)$ and $Y=\p(y)$ without constant term if and only if $X$ is
nilpotent in $A$.
\end{corollary}

\proof In one direction this is trivial, we should ensure   only
that if $X^N=0$ then $p(X,Y)^{2N}=0$ for any polynomial $p$ such
that $p(0,0)=0$ using the relation $XY-YX=Y^2$ which is an easy
check.$\Box$

\begin{theorem}\label{t2.1} Let $A_{\rho,n}$ be the image algebra of $R=\kk
xy\big/(xy-yx-y^2)$ under the $n$-dimensional representation
$\rho_n$ and $X=\rho_n(x)$, $Y=\rho_n(y)$ be its generators. Then
$A/J$ is a commutative one-generated ring $\k[x]/S(x)$, where
$S(x)=(x-\lambda_1)\dots(x-\lambda_k)$ and
$\lambda_1,\dots,\lambda_k$ are all different eigenvalues of the
matrix $X$. \end{theorem}

\proof From the corollary \ref{c2.2} we can see that $A/J$ is an
algebra of one generator $x$: $A/J \simeq \k[x]/I$. We are going to
find now an element which generates  the ideal $I$.

First of all by lemmas \ref{l2.2} and  \ref{l2.3} $S\in J(A)$, hence
$S(x)=(x-\lambda_1) \dots (x-\lambda_r) \in I$. Let we show now that
$S$ divides any element of $I$. If some polynomial $p\in \k[x]$ does
not vanish in some eigenvalue $\lambda$ of $X$ then $p(X)\notin
J(A)$. Indeed, the matrix $p(X)$ has a non-zero eigenvalue, than
$p(\lambda) \neq 0$ and hence $I-\frac 1{p(\lambda)}p(X)$ is
non-invertible. Therefore $p(X)\notin J(A)$. Thus, $S(x)$ is the
generator of $I$. This finishes the proof. $\Box$


\begin{corollary}\label{ct2.2} The system
$e_i=p_i(X)/p_i(\ll_i)$, where
$$p_i(X)=(X-\ll_1I)\dots\widehat{(X-\ll_iI)}\dots(X-\ll_rI)$$ and
$\ll_i$ are a different eigenvalues of $X=\p(x)$ is a complete
system of orthogonal idempotents of $A/J$.
\end{corollary}

\proof Orthogonality of $e_i$ is clear from the presentation of
$A/J$ as $\k[x]/{\rm id}(S)$ proven in theorem \ref{t2.1}. $\Box$

\begin{theorem}\label{t2.3}  For any finite dimensional representation
 $\p$ the semisimple part of $A_{\p}$ is a
product of a finite number of copies of the field $\k$:
$$ A/J = \prod_{i=1}^r \k_i, $$ where $r$ is the number of
different eigenvalues of the matrix $X=\rho(x)$.\end{theorem}

\proof We shall construct an isomorphism of $A/J$ and $\displaystyle
\prod_{i=1}^r \k_i$ using the system $e_i$, $i=1,\dots,r$ of
idempotents from the corollary  \ref{ct2.2}. Clearly $e_i$ form a
basis of $A/J$ as a linear space over $\k$. From the presentation of
$A/J$ as a quotient $\k[x]/{\rm id}(S)$ given in the theorem
\ref{t2.1} it is clear that the dimension of $A/J$ is equal to the
degree of polynomial $S(x)$, which coincides with the number of
different eigenvalues of the matrix $X$. Since idempotents $e_i$ are
orthogonal, they are linearly independent and therefore form a basis
of $A/J$. The multiplication of two arbitrary elements $a,b\in A/J$,
$a=a_1e_1+{\dots}+a_r e_r$, $b=b_1e_1+{\dots}+b_r e_r$ is given by
the formula $ab=a_1b_1e_1+{\dots}+a_r b_r e_r$ due to orthogonality
of the idempotents $e_i$. Hence the map $a\mapsto(a_1,\dots,a_r)$ is
the desired isomorphism of $A/J$ and $\displaystyle \prod_{i=1}^r
\k_i$.
$\Box$

Since all images turned out to be basic algebras we can associate to
them a {\it quiver} in a conventional way (see, for example,
\cite{Gr}, \cite{Kirich}).
The vertices will correspond to the idempotents $e_i$ or by the
corollary \ref{ct2.2} equivalently, to the different eigenvalues of
matrix $X$. The number of arrows from vertex $e_i$ to the vertex
$e_j$ is the ${\rm dim}_\k \, e_i(J/J^2)e_j$. There are a finite
number of such quivers in fixed dimension $n$ (the number of
vertices bounded by $n$, the number of arrows between any two
vertices roughly by $n^2$).

Let we prove now the following lemma. Denote by $\bar Y$ the image
of $Y$ under the factorization by the square of radical: $\bar Y
=\phi Y$, for $\phi:A \rightarrow A/J^2$.

\begin{lemma}\label{radrsq}
If in the representation $\p: R \rightarrow A, \, \,\, X=\p(x)$ has
only one eigenvalue $\lambda$, then the corresponding quiver $Q_A$
has one vertex and number of loops is a dimension of the vector
space ${\rm Sp}_k\{ \bar X-\lambda I, \bar Y\}$, which  does not
exceed 2.
\end{lemma}

\proof Due to the description of idempotents above in the case of
one eigenvalue the only idempotent is unit. Hence we have to
calculate ${\rm dim}_k J/J^2$, where $J=Jac(A)$. Since $X-\lambda I$
satisfy the same relation as $X$ we could apply the corollary
\ref{radft} and result immediately follows. $\Box$

 After we have proved that all image algebras are basic we can define an equivalence
 relation on
representations of RIT algebra using quivers of its images.

\begin{definition} Two representations $\rho_1$ and $\rho_2$ of the
algebra $R$ are quiver-equivalent $\rho_1 \sim_{Q} \rho_2$ if
quivers associated to algebras $\rho_1(R)$ and $\rho_2(R)$ coincide.
\end{definition}


As an example let us clarify the question on how many
quiver-equivalence classes appear in the family of representations
$$
{\cal M}_n=\{(X,Y)\in {\rm mod}(R,n)|{\rm rk}\,Y=n-1\}
$$
and which quivers are realized.

\begin{proposition}\label{t2.4}For any $n \geq 3$ families of representations ${\cal M}_n$ belong
to one quiver-equivalence class. Corresponding quiver consists of
one vertex and two loops. \end{proposition}

\proof This will directly follow from Lemma \ref{radrsq}, when we
ensure in section \ref{descript} that $X$ has only one eigenvalue in
the family ${\cal M}_n $ and take into account that when we have
full block $Y$, the dimension of the linear space ${\rm Sp}_k\{ \bar
X-\lambda I, \bar Y\}$ can not be smaller then 2. $\Box$

\section{Automorphisms of RIT algebras and multiplication formulas}

Here we intend to describe  the group of automorphisms of the RIT
algebra $R$ in order to use this information later on for the
classification results. It turned out to be quite small, compared
with automorphisms of the first Weyl algebra $A_1$, which contains
$R$ as a subalgebra. Automorphisms of the $A_1$  were described in
\cite{ML}, the case of an arbitrary Weyl algebra $A_n$ was discussed
in \cite{MK}.

First we shall prove lemmata on multiplication in RIT, it will be on
use for various purposes later on.

\subsection{Preliminary facts on multiplication in RIT}

Since the defining relation for R: $xy=yx+y^2$ form a Gr\"obner
basis with respect to the ordering $x>y$,
 the basis of our algebra as a vector space over $k$ consists of the
monomials $y^kx^l$, $k,l=0,1,\dots$. These are those monomials which
do not contain the highest term $xy$ of the defining relation.

We prefer to show this here in a canonical way. For this we shell
remaind the definition of a Gr\"obner  basis of an ideal and the
method of construction of a linear basis of an algebra given by
relations, based on the Gr\"obner basis technique. Using this
canonical method it could be easily shown that, for example, some
Sklyanin algebras enjoys a PBW property. This was proved in
\cite{od}, the arguments there are very intelligent and interesting
in their own right, but quite involved.

Let $A=k\langle X\rangle/I$. The first essential step is to fix an
ordering on the semigroup $\ss=\langle X \rangle$. We fix some
linear ordering in the set $X$.  Then we have to extend it to an
{\it admissible} ordering on $\ss$, i.e. it has to satisfy the
conditions:

1) if  $u,v,w \in \ss$ and $ u<v$ then $uw<vw$ and $wu<wv$

2) the descending chain condition (d.c.c.):
 there is no infinite properly
descending chain of elements of $\ss$.

We shall use the {\it deg\-ree-lex\-i\-c\-o\-g\-r\-a\-p\-h\-i\-cal}
ordering in the semigroup $\ss$, namely  for arbitrary $ \, u=
x_{i_1}\ldots x_{i_n}, v= x_{j_1}\ldots x_{j_k}\in\ss \quad \text{we
say} \quad
 u>v, \, \text{when either}\quad$ $deg\, u > deg\, v
$ $ \quad \text{or} \quad   deg\, u =  deg\, v \quad \text{and for
some} \,\, l: \,\, x_{i_l}>x_{j_l} \quad \text{and}\quad
x_{i_m}=x_{j_m}  \, \text{for any} \,\,  m<l.$
  This ordering is admissible.

Denote by $\bar f$ the highest term of polynomial $ f\in A=k\langle
X\rangle$ with respect to the above order.

{\bf Definition 4.2.} Subset $G \in I, I \triangleleft \langle X
\rangle$ is a {\it Gr\"obner basis} of an ideal if the set of
highest terms of elements of $G$ generates the ideal of highest
terms of $I: id\{\bar G\} = \bar I$.

{\bf Definition 4.3.} We will say that monomial $u \in \langle X
\rangle$ is {\it normal} if it does not contain as a submonomial any
highest term of an element of the ideal I.

From these two definitions it is clear that normal monomial is a
monomial which does not contain  any highest term of an element of
Gr\"obner basis of the ideal $I$. If Gr\"obner basis turned out to
be finite then the set of normal words is constructible.

In case when an ideal $I$ of defining relations for $A$ has a finite
Gr\"obner basis, the algebra called {\it standardly  finitely
presented(s.f.p.)}.

It is easy, but useful fact that $\langle X \rangle$ is isomorphic
to the direct sum $I \oplus {{\langle N \rangle}_{k}} $ as a linear
space over $k$, where $\langle N \rangle_{k}$ is the linear span of
the set of normal monomials from $\lxr$ with respect to the ideal
$I$. We claim  also (without check, which is not difficult) that the
set of normal words form a linear basis. Hence given a Gr\"obner
basis $G$ of an ideal $I$, we can construct a linear basis of an
algebra $A=\langle X \rangle / I$ as a set of normal (with respect
to $I$) monomials, at least in case when $A$ is s.f.p.

As a consequence we immediately get the following

\begin{lemma}\label{l4.1}  The system of monomials $y^n x^m$ form a basis of
algebra $R$ as a vector space over $k$.\end{lemma}

 We say that an element is in {\it normal form}, if
it is presented as a linear combination of normal monomials.

After we have a linear basis of normal monomials we should know how
to multiply them to get again an element in normal form.

Now we are going to prove  the following lemmata, where we express
precisely normal forms of some products.

\begin{lemma}\label{l3.1.1} The normal form of the monomial $xy^n$ in algebra $R$ is the following:

  $xy^n=y^nx+ny^{n+1}$.\end{lemma}

\proof This can be proven by induction on $n$. The case $n=1$ is
just our algebra's relation. Suppose $n>1$ and the equality
$xy^{n-1}=y^{n-1}x+(n-1)y^{n}$ holds. Multiplying it by $y$ from
the right and reducing by the relation $xy-yx=y^2$, we obtain
$$
xy^n=y^{n-1}xy+(n-1)y^{n+1}=y^nx+y^{n+1}+(n-1)y^{n+1}=y^nx+ny^{n+1}.
$$
The proof is now complete.$\Box$

\begin{lemma}\label{l3.1.2}
 The normal form of the monomial $x^ny$ in algebra $R$ is the following:
$x^ny=\sum\limits_{k=1}^{n+1}\alpha_{k,n}y^kx^{n-k+1}$, where
$\alpha_{k,n}= n! / (n-k+1)! $ for $ k=1,...,n+1$.
\end{lemma}

\proof  We are going to prove this formula inductively using the
previous lemma. As a matter of fact we shall obtain recurrent
formulas for $\alpha_{k,n}$. In the case $n=1$ the relation
$xy-yx=y^2$ implies the desired formula with
$\alpha_{1,1}=\alpha_{2,1}=1$. Suppose $n$ is a positive integer and
there exist positive integers $\alpha_{k,n}$, $k=1,\dots,n+1$ such
that $x^ny=\sum\limits_{k=1}^{n+1}\alpha_{k,n}y^kx^{n-k+1}$.
Multiplying the latter equality by $x$ from the left and using lemma
\ref{l3.1.1} we obtain
$$
x^{n+1}y=\sum_{k=1}^{n+1}\alpha_{k,n}xy^kx^{n-k+1}=
\sum_{k=1}^{n+1}\alpha_{k,n}y^kx^{n-k+2}+
\sum_{k=1}^{n+1}\alpha_{k,n}ky^{k+1}x^{n-k+1}.
$$
Rewriting the second term as
$\sum\limits_{k=1}^{n+2}\alpha_{k-1,n}(k-1)y^{k}x^{n-k+2}$ (here
we assume that $\alpha_{0,n}=0$), we arrive to
$$
x^{n+1}y=\sum_{k=1}^{n+2}\alpha_{k,n+1}y^kx^{n-k+2},
$$
where $\alpha_{k,n+1}=\alpha_{k,n}+(k-1)\alpha_{k-1,n}$ for
$k=1,\dots,n+1$ and $\alpha_{n+2,n+1}=(n+1)\alpha_{n+1,n}$.

Let we prove now the formula for $\alpha_{k,n}$. For $n=1$ it is
true since $\alpha_{1,1}=\alpha_{1,2}=1$. Then we use  inductive
argument. Suppose the formula is true for $n$. We are going to apply
the recurrent formula appeared above:
$$ \alpha_{k,n+1}=\alpha_{k,n}+(k-1)\alpha_{k-1,n}=
\frac{n!}{(n-k+1)!}+(k-1)\frac{n!}{(n-k+2)!}=
\frac{(n+1)!}{(n-k+2)!} $$ and the formula is verified for $1\leq
k\leq n+1$.  For $k=n+2$, we have
$\alpha_{n+2,n+1}=(n+1)\alpha_{n+1,n}=(n+1)n!=(n+1)!$. This
completes the proof. $\Box$

\subsection{Automorphisms of RIT algebras }\label{aut}

 We  are going to describe the automorphism group of the simplest RIT algebra here.
 We shall prove.

\begin{theorem}\label{taut}
 All automorphisms of $R=k\langle x,y|xy-yx=y^2\rangle$
are of the form $x\mapsto \alpha x+p(y)$, $y\mapsto \alpha y$, where
$\alpha\in k\setminus\{0\}$ and $p\in k[y]$ is a polynomial on $y$.
Hence the group of automorphisms isomorphic to a semidirect product
of an additive group of polynomials $k[y]$ and a multiplicative
group of the field  $k^*: \,\,{\rm Aut}(R) \simeq k[y]
\leftthreetimes k^*$.
\end{theorem}

\proof Key observation for this proof is that in our algebra there
exists a minimal ideal with commutative quotient. Namely, the
two-sided ideal $J$ generated by $y^2$.

\begin{lemma}\label{lcommq} If the quotient $R/I$ is commutative then $y^2\in I$ (that
is $J\subseteq I$). \end{lemma}

\proof The images of $x$ and $y$ in this quotient commute. Hence
$$
0=(x+I)(y+I)-(y+I)(x+I)=xy-yx+I=y^2+I.
$$
Therefore $y^2\in I$. $\Box$

The property of an ideal to be a minimal ideal with commutative
quotient is invariant under automorphisms.

Let us denote by $\widetilde y=f(x,y)$ the image of $y$ under an
automorphism $\phi$. Then the ideal generated by $\widetilde y^2$
coincides with the ideal generated by $y^2$: $J=\langle
y^2\rangle=\langle \widetilde y^2\rangle$.

Using the  property of multiplication in $R$ from lemma
\ref{l3.1.2}, we can see that two-sided ideal generated by $y^2$
coincides with the left ideal generated by $y^2$: $Ry^2R=y^2R$.
Indeed, let us present an arbitrary element of $Ry^2R$ in the form
$\sum a_iy^2b_i$, where $a_i$, $b_i\in R$ are written in the normal
form $a_i=\sum \alpha_{k,l}y^kx^l$, $b_i=\sum \beta_{k,l}y^kx^l$.
Using the relations from lemma \ref{l3.1.1}, we can pull $y^2$ to
the left through $a_i$'s and get the sum of monomials, which all
contain $y^2$ at the left hand side. Thus,  $\sum a_iy^2b_i=y^2u$,
$u\in R$.

Obviously automorphism maps the one-sided ideal $y^2R$ onto the
one-sided ideal $\widetilde y^2R$, both of which coincide with
$J=\langle y^2\rangle=\langle \widetilde y^2\rangle$. From this we
obtain a presentation of $y^2$ as $\widetilde y^2 u$ for some $u\in
R$. Considering  usual degrees of these polynomials (on the set of
variables $x,y$), we get $2=2k+l$, where $k={\rm deg}\,\widetilde y$
and $l={\rm deg}\,u$. Obviously $k\neq 0$. Hence the only
possibility is $k=1$ and $l=0$.

Thus, $\phi(y)=\widetilde y=\alpha x+\beta y+\gamma$ and $u=c$ for
some $\alpha,\beta,\gamma,c\in k$. Substituting these expressions
into the equality $y^2=\widetilde y^2 u$, we get $c(\alpha x+\beta
y+\gamma)^2=y^2$. Comparing the coefficients of the normal forms
of the right and left hand sides of this equality, we obtain
$\alpha=\gamma=0$, $\beta\neq 0$. Hence $\phi(y)=\beta y$.

Now we intend to use invertibility of $\phi$. Due to it there
exists $\alpha_{ij}\in k$ such that $x=\sum \alpha_{ij}\widetilde
y^i \widetilde x^j$. Substituting $\widetilde y=\beta y$, we get
$\ddd x=\sum_{r=0}^N p_r(y)\widetilde x^r$, where $N$ is a
positive integer, $p_r\in k[y]$ and $p_N\neq 0$. Comparing the
degrees on $x$ of the left and right hand sides of the last
equality we obtain $1=kN$, where $k={\rm deg}_x\widetilde x$.
Hence $k=N=1$, that is $x=p_0(y)+p_1(y)\widetilde x$ and
$\widetilde x=q_0(y)+q_1(y)x$, where $p_0,p_1,q_0,q_1\in k[y]$.
Substituting $\widetilde x=q_0(y)+q_1(y)x$ into
$x=p_0(y)+p_1(y)\widetilde x$, we obtain $q_1\in k$, that is
$\widetilde x=cx+p(y)$ for $c\in k$. One can easily verify that
the relation $\widetilde x\widetilde y-\widetilde y\widetilde x=
\widetilde y^2$ is satisfied for $\widetilde x=cx+p(y)$,
$\widetilde y=\beta y$ if and only if $c=\beta$. This gives us the
general form of the automorphisms:  $\widetilde x=cx+p(y)$,
$\widetilde y=c y, c\neq 0$.

Now we see that the group of automorphisms is a semidirect product
of the normal subgroup isomorphic to the additive group of
polynomials $k[y]$ and the subgroup isomorphic to the multiplicative
group $k^*$. The precisely written formula for multiplication in
${\rm Aut} R$ is the following:

$$\phi_1 \phi_2 = (p_1(y),c_1) (p_2(y), c_2) = (c_2
p_1(y)+p_2(c_1y), c_1c_2)$$

for $\phi_1, \phi_2 \in {\rm Aut} R $. $\Box$

\section{Irreducible modules,
description of all finite dimensional modules}\label{descript}

We intend to prove here the following

\begin{theorem}\label{tdm} The description of the complete set of finite
dimensional representations of $R$ (subject to the Jordan form of
$Y$) are given by

{\small
\begin{equation}
\label{rho}
Y_n=\!\left(\Mfour\right)\!,\ X_n=\!\left(\Meight\right)
\end{equation}}
where partition on blocks in  $X_n$ correspond to the partition
defined by the Jordan form of $Y_n$. Diagonal blocks of $X_n$ are
matrices $X_n^0+T$, where $X^0$ is a matrix with the vector
$[0,1,2,...]$ on the first upper diagonal and zeros elsewhere, $T$
is an arbitrary upper diagonal Toeplitz matrix. All the rest of
blocks of $X$ are upper diagonal rectangular Toeplitz matrices.

\end{theorem}

From this theorem immediately follows a precise description of all
{\it irreducible} and {\it completely reducible} modules.

\begin{corollary}\label{c3.1}
A complete set of pairwise  non-isomorphic finite dimensional {\it
irreducible} $R$-modules is $\{ S_{a} | {a} \in k\}$, where $S_{a}$
defined by the following action of $X$ and $Y$ on one-dimensional
vector space: $Xu=\a u, Yu=0.$


All {\it completely reducible} representations are given by
matrices: $Y_n=(0)$, $X_n$ is a diagonal matrix ${\rm diag}
(a_1,...,a_n)$.
\end{corollary}


\proof Let we describe an arbitrary representation $\rho_n:R\to
M_n(k)$ of $R$, for $n\in\N$. We can assume that the image of one of
the generators $Y=\rho_n(y)$ is in normal Jordan form.

{\it Full Jordan block case.}

 Let us first find all possible matrices
$X=\rho_n(x)$ in the case when $Y$ is just full Jordan block:
$Y=J_n$. We have to find than matrices $X=(a_{ij})$ satisfying the
relation $[X,J_n]=J_n^2$. Let $B=[X,J_n]=(b_{ij})$, then
$b_{ij}=a_{i+1,j}-a_{i,j-1}$. From the condition $B=J_n^2$ it
follows that $b_{ij}=0$ if $i\neq j-2$ and $b_{ij}=1$ if $i=j-2$.
Here and later on we will  use the following numeration of
diagonals: main diagonal has number 0, upper diagonals have positive
numbers $1,2,\dots,n-1$ and lower diagonals have negative numbers
$-1,-2,\dots,-n+1$:

$$
\left(\Mone\right)
$$

\medskip

The first condition above means than that in the matrix $X$ elements
of any diagonal with number $0 \leq k\neq 1$ coincide and are zero
for $k<0$. From the second condition it follows that the elements of
the first upper diagonal form an arithmetic progression with
difference 1:  $a+1,\dots,a+n-1$.

Therefore we have the following sequence of
representations:

\begin{equation}
\label{eps}
Y_n=\left(\MtwoA\right),\quad X_n=\left(\MtwoB\right)
\end{equation}

Here and below we will draw a diagonal as a continuous line if
all its elements coincide and as a thick line if its elements form
am arithmetic progression with difference one.

Denote by $X^0$ a matrix with the sequence $0,1,2,...$ on the first
diagonal and zeros elsewhere.  Then our family of representations
consists of pairs of matrices $(X_n,Y_n)=(X_n^0+T, J_n)$, where $T$
is an arbitrary upper diagonal Toeplitz matrix. Let we remind that
{\it upper diagonal (rectangular) Toeplitz matrix} is a matrix with
entries $a_{ij}$ defined only by the difference $i-j$. It has zeros
below the main diagonal (or upper main diagonal in a proper
rectangular case).


Note that one could get a clue on what the set of representations is
from the following observation. First, the matrix $X^0$ satisfies
the relation $[X^0,Y]=Y^2$ for $Y=J_n$. On the other hand a matrix
$X=X^0+M$ satisfies the relation $[X,Y]=Y^2$ if and only if $M$
commutes with $Y=J_n$. Any matrix having only one non-zero diagonal
with equal elements on it commutes with $Y=J_n$. Hence we have at
least all linear combinations of those matrices in the set of
representations, additional arguments as above show that there are
no others.

{\it The case of an arbitrary partition}.

 Consider now the
general case when the Jordan normal form of $Y$ contains several
Jordan blocks: $Y=(J_1,...,J_m)$.


Cut an arbitrary matrix $X$ into the square and rectangular blocks
of corresponding size, denote the blocks by $A_{ij}, i,j={\overline
1,m}$.


Then we can describe the structure of the matrix $B=[X,Y]$ in the
following way:

{\small

$$
B=\left(\Msix\right),\ \ \text{where\ }[A_{ij}]=A_{ij}J_i-
J_jA_{ij}.
$$}
From the condition $B=Y^2$ we have that $[A_{ii},J_i]=J_i^2$ and
hence $A_{ii}$ is the same as in the previous case when $Y$ was just
a full Jordan block and $A_{ij}J_i-J_jA_{ij}=0$ for $i\neq j$. The
latter condition means that $A_{ij}$ for $i\neq j$ has the following
structure.

{\small
$$
\left(\MsevenA\right)\quad\text{or}\quad
\left(\MsevenB\right)
$$}

 The elements of any diagonal here marked as a line are equal to
each other, elements of different diagonals could be different and
they are equal to zero below the upper diagonal of maximal length
(the matrix is non-square in general). As a result we have the
family of representations described in the theorem \ref{tdm}. $\Box$



\section{$R$ is residually finite dimensional
}\label{rfd}

Let us consider now one of the sequences of representations
constructed in the  previous section: $\ee_n: R \rightarrow {\rm
End}\,\, k^n$, defined by $\ee_n(y)=J_n$, $\ee_n(x)=X_n^0$ as above.
Note that this sequence is basic in the following sense. As was
actually shown in \ref{descript}, all representations (\ref{eps})
corresponding to $Y$ with one Jordan block could be obtained from
$\ee_n$ by the following automorphism of $R$, \, $\phi:R
\longrightarrow R: x \mapsto x+a, y \mapsto y$ where $a \in R$ such
that $[a,y]=0$.





In addition to the conventional equivalence relation on the
representations given by simultaneous conjugation of matrices:
$\rho' \sim \rho''$ if there exists $g \in GL(n)$ such that
$g\rho'g^{-1}=\rho''$ or equivalently, R-modules corresponding to
$\rho'$ and $\rho''$ are isomorphic, we introduce here one more
equivalence relation.

\begin{definition} We say that two representations of the
algebra $R$ are {\it auto-equivalent} (equivalent up to
automorphism) $\p' \sim_A \p'' $ if  there exists $\phi \in \rm{Aut}
(R)$ such that $ \p' \phi \sim \p''$.
\end{definition}

So  we can state that any full block representation
 is auto-equivalent to $\ee_n$ for appropriate $n$.

We will prove now that the sequence of representations $\ee_n$
asymptotically is faithful.

Start with the calculation of matrices which are image of monomials
$y^k x^m$ under representation $\ee_n$.

\begin{lemma}\label{lImyx} For the representation $\ee$ as above
the matrix $\ee(y^k x^m)$ has the following shape: there is only one
nonzero diagonal, number $k+m$, in the above numeration, where
appears the sequence $p(0),p(1),...,p(j),...$ of values of degree
$m$ polynomial $p(j)=(k+j)...(k+m+j-1)=\prod_{i=1}^m (k+j+i)$.
\end{lemma}

\proof Image $\ee(x^m)$ of the monomial $x^m$ is a matrix with
vector $[1 \cdot 2 \cdot...\cdot m, 2 \cdot 3 \cdot ...\cdot (m+1),
...] $ on the (upper) diagonal number $m$ in the above numeration
and zeros elsewhere. Multiplication by  $\ee (y^k)$ acts on matrix
by  moving up all rows on $k$ steps. We can  now see that matrix
corresponding to the polynomial $y^k x^m$ can have only one nonzero
diagonal, number $m+k$, and vector in this diagonal is the
following: $[(k+1)...(m+k), (k+2)...(m+k+1),...]$.$\Box$

\begin{theorem}\label{t4.1} Let $\ee_n$ be the sequence of representations of $R$ as above.
Then $\cap_{n=0}^\infty \, {\rm ker} \, \ee_n = 0$. \end{theorem}

\proof We are going to show that  $\ee_n(f) \neq 0$ for $n \geq 2
\deg f$. Suppose that $n$ is sufficiently large and  $\ee_n(f)$ is
zero and get a contradiction. Denote by $l$  degree of polynomial
$f$, and let $f=f_1 +...+f_l$ be a decomposition of $f \in R$ on the
homogeneous components of degrees $i=1,...,l$ respectively. From
lemma \ref{lImyx} we know now how the matrix which is an image of an
arbitrary monomial $y^k x^m$ looks like.

Applying the lemma \ref{lImyx} to each homogeneous part of the given
polynomial $f$ we get
$$
f_l=\sum_{k+m=l}a_{k,m}y^kx^m=\sum_{r=0}^la_ry^{l-r}x^r
$$
is a sum of matrices $\displaystyle \sum_{r=0}^la_rM_r$, where $M_r$
has the vector $[(p(0),...,p(j)]$:
$$
\left(\prod_{i=1}^r(l-r+i),\prod_{i=1}^r(l-r+i+1),\dots,
\prod_{i=1}^r(l-r+i+j),\dots \right)
$$
on the diagonal number l (all other entries are zero). The number on
the $j$-th place of this diagonal is the value in $j$ of the
polynomial
$$P(j)=
(l-r+j) \cdot ...\cdot (l+j-1)$$ of degree exactly $r$. Therefore
the sum $\displaystyle \sum_{r=0}^l a_rM_r$ has a polynomial on $j$
of degree $N=\max\{r:a_r\neq0\}$ on the diagonal number $l$. Since
any polynomial of degree $N$ has at most $N$ zeros we arrive to a
contradiction in the case when $l$th diagonal has length more than
$l$. Hence for any $n \geq 2 \deg f$, \, $\ee_n(f) \neq 0$. $\Box$

Let we recall that an algebra $R$ {\it residually has some property
$\PP$} means that there exists a system of equivalence relations
$\tau_i$ on $R$ with trivial intersection, such that in the quotient
of $R$ by any $\tau_i$ property $\PP$ holds.

From the Theorem \ref{t4.1}  we have the following corollary
considering equivalence relations modulo ideals ker $\ee_n$.

\begin{corollary}\label{c4.1} Algebra $R$ is residually finite dimensional.\end{corollary}

\section{Indecomposable modules}\label{ind}

\begin{lemma}\label{lam}
Let $M=(X,Y)$ be a (finite dimensional) indecomposable module over
$R=k\langle x,y|xy-yx=y^2\rangle$. Then $X$ has a unique eigenvalue.
\end{lemma}

\proof Denote by $M_\lambda^X$ the main eigenspace for $X$
corresponding to its eigenvalue $\lambda$:
$M_\lambda^X=\bigcup\limits_{k=0}^\infty {\rm ker}\,(X-\lambda
I)^k$. Obviously $M_\lambda^X={\rm ker}\,(X-\lambda I)^m$, where $m$
is the maximal size of blocks in the Jordan normal form of $X$. It
is well-known that $M=\mathop{\oplus}\limits_{i} M_{\lambda_i}^X$,
where the direct sum is taken over all different eigenvalues
$\lambda_i$ of $X$. We shall show that $M_{\lambda_i}^X$ are in fact
$R$-submodules.

Let $u\in M_{\lambda}^X$, that is $(X-\lambda I)^mu=0$. We calculate
$(X-\lambda I)^nYu$ for arbitrary $n$. Using the fact that the
mapping defined on generators $\phi(x)=x-\lambda$, $\phi(y)=y$
extends to an automorphism of $R$ (see \ref{aut}), we can apply it
to the multiplication  formula from Lemma \ref{l3.1.2} to get
$(x-\lambda)^ny=\sum\limits_{k=1}^{n+1}y^k(x-\lambda)^{n-k+1}$.
Taking into account that $Y^l=0$ for some positive integer $l$, we
can choose $N$ big enough, for example  $N\geq m+l$, such that


$$
(X-\lambda I)^NYu=\sum_{k=1}^{N+1}\alpha_{k,N}Y^k(X-\lambda
I)^{N-k+1}u=0
$$

either due to $(X-\lambda I)^{N-k+1}u=0$ or due to $Y^k=0$.

 This shows that $Yu\in M_\lambda^X$, that is
$M_\lambda^X$ is invariant with respect to $Y$. $\Box$

As an immediate corollary we have the following.

\begin{proposition}\label{pdecomp} Any finite dimensional $R$-module $M$ decomposes into the
direct sum of submodules  $M_{\lambda_i}^X$ corresponding to
different eigenvalues $\lambda_i$ of $X$. \end{proposition}


\begin{corollary}\label{cindlocal} Let $M$ be indecomposable module corresponding to the
representation $\p: R \longrightarrow {\rm End}(k^n)$, and $A_n$ is
the image of this representation. Than $A_n$ is local algebra, e.i.
$A_n / J(A_n) = k$. \end{corollary}

\proof This follows from the above lemma \ref{lam} and fact that any
image algebra is basic  with semisimple part  isomorphic to the sun
of $r$ copies  of the field $k$: $\oplus_{r} k $, where $r$ is a
number of different eigenvalues of $X$, which was proved in the
\ref{t2.3}.$\Box$

Now using the definition of quiver for the image algebra given in
section \ref{quiv}  and lemma \ref{radrsq} we give a complete
description of quiver equivalence classes of indecomposable modules.

\begin{corollary}\label{cindQuiv} Quiver corresponding to the indecomposable module has one
vertex. The number of loops is one or two, which is a dimension of
the vector space ${\rm Sp}_k\{ \bar X-\lambda I, \bar Y\}$, where
$\bar  X=\phi X$, $\bar  Y=\phi Y$ for $\phi: A \rightarrow A/J^2$.
\end{corollary}

As another consequence of the proposition \ref{pdecomp} we can
derive an important information on how  to glue irreducible modules
to get indecomposables. It turned out that it is possible to glue
together nontrivially  only the copies of the same irreducible
module $S_a$.


\begin{corollary}\label{ext} For arbitrary non-isomorphic irreducible   modules
$S_a, S_b$,  $${\rm Ext}^1_k (S_a,S_b)=0, \,\,{\rm if}\,\, a \neq
b.$$
\end{corollary}
\proof Indeed, in corollary \ref{c3.1} we derive that irreducible
module $S_i$ is one dimensional and given by $X=(a), Y=(0)$, $a \in
k$. If $a \neq b$ then for $[M] \in {\rm Ext}^1_k (S_a,S_b)$,
corresponding $X$ has two different eigenvalues, namely $a$ and $b$.
Then by the above lemma $M$ is decomposable and $[M]=0$. $\Box$

\section{Equivalence of some subcategories in mod$\,R$}\label{eqc}

Let we denote by mod $R(\lambda)$ the full subcategory in mod$\,R$
consisting of modules with the  unique eigenvalue $\lambda$ of $X$:
${\rm mod}\,R(\lambda)=\{M\in {\rm mod}\,R|M=M_\lambda(X)\}$. Let us
define the functor $F_\lambda$ on ${\rm mod}\,R$, which maps a
module $M$ to the module $M_\lambda$ with the following new action
$rm=\phi_\lambda(r)m$, where $\phi_\lambda$ is an automorphism of
$R$ defined by $\phi_\lambda(x)=x+\lambda$, $\phi_\lambda(y)=y$. The
restriction of $F_\lambda$ to ${\rm mod}\,R(\lambda) $ is an
equivalence of categories $F_\lambda:{\rm mod}\,R(\lambda)\to {\rm
mod}\,R(\mu+\lambda)$ for any $\mu\in k$. In particular, we have an
equivalence of the categories ${\rm mod}\,R(\lambda)$ and ${\rm
mod}\,R(0)$.

To use this equivalence of categories it is necessary to know that
in most cases (but not in all of them), the eigenvalues of the
matrix $X$ are just entries of the main diagonal in the standard
shape of the matrix described in the Theorem \ref{tdm}, more
precisely.

\begin{theorem}\label{tnn} Let in the
basis {\it E} of the representation vector space, $Y$ is in the
Jordan normal form, and Jordan blocks have pairwise different sizes:
$n_1,n_2,\dots,n_k$. Then in the same basis $X$ has the shape
(\ref{rho}) with  numbers $\lambda_1,\dots,\lambda_k$ on the
diagonals of the main blocks, where $\lambda_j$ are eigenvalues of
$X$ (not necessarily different).
\end{theorem}

\proof Let we first introduce the denotation for the basis $\cal E$:

$$
e^{1,1},\dots,e^{1,n_1},e^{2,1},\dots,e^{2,n_2},\dots,
e^{k,1},\dots,e^{k,n_k}.
$$

Consider the set $\cal A$ of the matrices (in the same basis) such
that $A_{(j,l),(j,l)}=c_j$, $1\leq j\leq k$, $1\leq l\leq n_j$ and
$A_{(i,s),(j,l)}=0$ if $n_j<n_i$ and $l>s-n_j$ and if $n_j>n_i$ and
$l>s$. One can easily verify that $\cal A$ is an algebra with
respect to the matrix multiplication. Let also $\cal D$ be the
subalgebra of diagonal matrices in ${\cal A}$ and $\phi:{\cal A}\to
{\cal D}$ be the natural projection ($\phi$ acts annihilating the
off-diagonal part of a matrix).

Looking at the multiplication in $\cal A$ it is straightforward to
see that $\phi$ is an algebra morphism, that is $\phi(I)=I$,
$\phi(AB)=\phi(A)\phi(B)$ and $\phi(A+B)=\phi(A)+\phi(B)$. It is
also easy to check, calculating the powers of the matrix, that if
$A\in\cal A$ and $\phi(A)=0$ then the matrix $A$ is nilpotent. Since
the matrices of the form (\ref{rho}) belong to $\cal A$, it suffices
to verify that the eigenvalues of any $A\in \cal A$ coincide with
the eigenvalues of $\phi(A)$.

First, suppose that $\lambda$ is not an eigenvalue of $A$. That is
the matrix $A-\lambda I$ is invertible: there exists a matrix
$B\in\cal A$ such that $(A-\lambda I)B=I$. Here we use the fact that
if a matrix from a subalgebra of the matrix algebra is invertible,
then the inverse belongs to the subalgebra. Then $\phi((A-\lambda
I))\phi(B)=\phi((A-\lambda I)B)=\phi(I)=I$. Therefore $\lambda$ is
not an eigenvalue of $\phi(A)$. On the other hand, suppose that
$\lambda$ is not an eigenvalue of $\phi(A)$. Then $\phi(A)-\lambda
I$ is invertible. Clearly
$$
A-\lambda I=(\phi(A)-\lambda I)(I+(\phi(A)-\lambda
I)^{-1}(A-\phi(A))).
$$
Let $B=(\phi(A)-\lambda I)^{-1}(A-\phi(A))$. Since $\phi$ is a
projection, we have that
$$
\phi(B)=\phi((\phi(A)-\lambda I)^{-1})(\phi(A)-\phi(A))=0.
$$
As we have already mentioned this means that the matrix $B$ is
nilpotent and therefore $I+B$ is invertible. Hence $A-\lambda
I=(\phi(A)-\lambda I)(I+B)$ is invertible as a product of two
invertible matrices. Therefore $\lambda$ is not an eigenvalue of
$A$. Thus, eigenvalues of $A$ and $\phi(A)$ coincide. This completes
the proof. $\Box$

\section{Analogue of the Gerstenhaber theorem for commuting
matrices}\label{Ger}

In this section we intend to prove an analog of the
Gerstenhaber-Taussky-Motzkin theorem (see \cite{Ger}, \cite{MT},
\cite{Gur}) on the dimension of images of representations of two
generated algebra of commutative polynomials $k[x,y]$. This theorem
says that any algebra generated by two matrices $A,B\in M_n(k)$ of
size $n$ which commute $AB=BA$ has dimension not exceeding $n$. We
consider instead of commutativity the relation $XY-YX=Y^2$ and prove
the following

\begin{theorem}\label{tG}
 Let $\rho_n:R\to M_n(k)$ be an arbitrary
$n$-dimensional representation of $R=k \langle x,y|xy-yx=x^2
\rangle$ and $A_n=\rho_n(R)$ be the image algebra. Then the
dimension of $A_n$ does not exceed $\frac{n(n+2)}4$ for even $n$ and
$\frac{(n+1)^2}4$ for odd $n$.

This estimate  is optimal and attained for the image algebra
corresponding to full-block $Y$.

\end{theorem}



We divide the proof in two lemmas. Start with the second statement
of the theorem, that is calculation of the dimension of image
algebras in full-block case.

\begin{lemma}\label{tGfull-bl}
 Let $X,Y\in M_n(k)$ be matrices of the size $n$ over the field
$k$, satisfying the relation $XY-YX=Y^2$  and $Y$ has as a Jordan
normal form one full block. Denote by ${A}$ the algebra generated by
$X$ and $Y$. Then for odd $n=2m+1$, ${\rm dim}\,
A=\frac{(n+1)^2}{4}$ and for even $n=2m$, ${\rm dim}\,
A=\frac{n(n+2)}{4}$.

\end{lemma}

\proof In the Lemma \ref{lImyx} we already have computed the
matrices, which are images of monomials $y^kx^m$ under the
representation $\ee: (x,y)\mapsto (X^0,J_n)$. Due to the fact that
any representation $\p$ with full block $Y$ could be obtained from
$\ee$ by composition with the R-automorphism $\phi: x\mapsto x+a,
y\mapsto y, $  where $[a,y]=0$, it is enough to calculate the
dimensions of images for $\ee$.

Let we recall how matrices $\epsilon(y^kx^m)$ look like and
calculate here the dimension of their linear span.

 The matrix
$\epsilon(y^{l-r}x^r)$ on the $l$-th upper diagonal has a vector
$(p(0),p(1),\dots)$, where
$$
p(j)=(l-r+j)\dots(l+j-1)=\prod_{i=1}^r(l+j-r+i)
$$
and zeros elsewhere. In the $j$-th place of the $l$-th diagonal we
have a value of a polynomial of degree exactly $r$. Those diagonals
which have number less then  the number of elements in it give the
impact to the dimension equal to the dimension of the space of
polynomials of corresponding degree. When the diagonals become
shorter (the number of elements less then the number of the
diagonal) then the impact to the dimension of this diagonal equals
to the number of the elements in it. Thus, if $n=2m+1$, ${\rm dim}\,
A=1+\dots+m+(m+1)+m+\dots+1=(m+1)^2=\frac{(n+1)^2}{4}$. When $n=2m$,
we have ${\rm dim}\, A=1+\dots+m+m+\dots+1=m(m+1)=\frac{n(n+1)}{4}$.
$\Box$

\subsection{Maximality of the dimension in the full-block case}

\def\frac#1#2{{#1}\over{#2}}
\def\text#1{{\rm {#1}}}

We know now that any representation of $R$,
which is isomorphic to one with $Y$ in full-block Jordan normal form
gives us as an image  the same algebra described in lemma
\ref{tGfull-bl} as a certain set of matrices, of dimension
${n(n+2)\over 4}$ for even $n$ and $\frac{(n+1)^2}4$ for odd $n$.
We intend to prove that this dimension is maximal among dimensions
of all image algebras for arbitrary representation, that is the
first part of the theorem \ref{tG}.

We start with the proof  that this dimension is an upper bound for
any image algebra of indecomposable representation.

The simple preliminary fact we will need is the following.

\begin{lemma}\label{ltr}
Matrices  $X$ and $Y$ satisfying the relation $XY-YX=Y^2$ can be by
simultaneous conjugation brought to a triangular form.
\end{lemma}

\proof
Using the defining relation and the fact that $Y$ is nilpotent we
can see that any eigenspace of $Y$ is invariant under $X$. Hence $X$
and $Y$ has joint eigenvector $v$. Then we consider quotient
representation on the space $V/\{v\}$ which has the same property.
Continuation of this process supply us with the basis where both $X$
and $Y$ are triangular.  $\Box$

\begin{lemma}\label{tGerst}
Let $\rho_n:R\to M_n(k)$ be an indecomposable  $n$-dimensional
representation of $R$ and $A_n=\rho_n(R)$ be the image algebra. Then
the dimension of $A_n$ does not exceed $\frac{n(n+2)}4$ for even $n$
and $\frac{(n+1)^2}4$ for odd $n$. \rm
\end{lemma}

\proof The algebra $A_n=\{\sum \alpha_{k,m}Y^kX^m\}$ consists now of
triangular matrices. Let we present the linear space $UT_n$ of upper
triangular $n\times n$ matrices as the direct sum of two subspaces
$UT_n=L_1\oplus L_2$, where $L_1$ consists of matrices with zeros on
upper diagonals with numbers $l,\dots,n$ and $L_2$ consists of
matrices with zeros on upper diagonals with numbers $1,\dots,l-1$,
where $l=(n+1)/2$ for odd $n$ and $l=n/2+1$ for even $n$. Let $P_j$,
$j=1,2$ be the linear projection in $UT_n$ onto $L_j$ along
$L_{3-j}$. Since $A_n$ is a linear subspace of $UT_n$, we have that
$A_n\subset M_1+M_2$, where $M_j=P_j(A_n)$. Therefore $\text{dim}\,
A_n\leq \text{dim}\,M_1+\text{dim}\,M_2$. The dimension of $M_1$
clearly does not exceed the dimension of the linear span of those
matrices $Y^kX^m$, which do not belong to $L_2$. Thus,
$$
\text{dim}\, M_1\leq \text{dim}\, \langle Y^kX^m| k+m<l-1\rangle _k.
$$

Here we suppose that $X$ (as well as $Y$) is nilpotent. We can do
this because the module is indecomposable. Indeed, the lemma
\ref{lam} says that for an indecomposable module $X$ has a unique
eigenvalue. This implies that any indecomposable representation is
autoequivalent to one with nilpotent $X$ and $Y$ due to the
automorphism of $R$ defined by $\phi_{\lambda}(x)=x-\lambda$,
$\phi_{\lambda}(y)=y$. Since autoequivalent representations has the
same image algebras we can suppose that $X$ is nilpotent.

Thus the dimension of $ M_1$ does not exceed the number of the pairs
$(k,m)$ of non-negative integers such that $k+m<l-1$, which is equal
to $1+\dots+(l-1)$. On the other hand $\text{dim}\,M_2\leq
\text{dim}\,L_2$ and the dimension of $L_2$ does not exceed the
total number of entries in the non-zero diagonals.
$$
{\rm dim}\,M_2\leq 1+\dots + (n-l+1).
$$
Taking into account that $\text{dim}\, A_n\leq
\text{dim}\,M_1+\text{dim}\,M_2$, we have
$$
{\rm dim}\, A_n\leq 1+\dots+(l-1)+1+\dots+(n-l+1).
$$
The latter sum equals $\frac{n(n+2)}4$ for even $n$ and
$\frac{(n+1)^2}4$ for odd $n$.
$\Box$

After we have proved the estimation for the indecomposable modules,
it is easy to see that the same estimate holds for arbitrary module,
since the function $n^2$ is convex.

On the other hand as it was shown in the Lemma \ref{tGfull-bl} this
estimate is attained on the algebra $A_n=\epsilon(R)$ in the case of
a full-block $Y$.
This completes the proof of the theorem \ref{tG}.



\section{
Parametrizable families of representations}\label{param}

Here we suppose that $k=\C$.
 Let we consider the variety of
$R$-module structures on $k^n$ and denote it by $mod(R,n)$. Such
structures are in 1-1 correspondence to a $k$-algebra homomorphisms
$R\to M_n(k)$ ($n$-dimensional representations), or equivalently to
a pairs of matrices $(X,Y)$, $X,Y\in M_n(k)$, satisfying the
relation $XY-YX=Y^2$. The group $GL_n(k)$ acts on $mod(R,n)$ by
simultaneous conjugation and orbits of this action are exactly the
isomorphism classes of $n$-dimensional $R$-modules. Denote this
orbit of a module $M$ or of a pair of matrices $(X,Y)$ as ${\cal
O}(M)$ or ${\cal O}(X,Y)$ respectively. Consider also the following
stratification. Let ${\cal U_P}$ be the set of all pairs $(X,Y)$
satisfying the relation, where $Y$ has a fixed Jordan form. Here
$\cal P$ stands for the partition of $n$, which defines the Jordan
form of $Y$. Clearly ${\cal U_P}$ is a union of all orbits where $Y$
has a Jordan form defined by partition $\cal P$: $ {\cal
U_P}=\bigcup\limits_{Y{\rm with  \, Jordan \, form } \atop {\rm
defined \, by \, the \, partition \, {\cal P}} } {\cal O}(X,Y). $ We
will write ${\cal U}_{(n)}$ for the stratum corresponding to $Y$
with the full Jordan block: ${\cal P}=(n)$.

Another action involved here is an action of the subgroup of $GL_n$
on those pairs $(X,Y)$, where $Y=J_{\cal P}$ is in fixed Jordan
form. Denote this space by $W_{\cal P}$. The subgroup which acts
there is clearly the centralizer of the given Jordan matrix:
$Z(J_{\cal P})$. Orbits of the action of $Z(J_{\cal P})$ on the
space $W_{\cal P}$ are just parts of orbits above: ${\cal
O_P}(X)={\cal O}(X,Y) \cap W_{\cal P}$.

We suggest here to consider in stead of action of $GL_n$ on the
whole space an action of centralizer $Z(J_{\cal P})$ on the smaller
space $W_{\cal P}$. While the group which acts is not reductive any
more and has a big unipotent part, we act just on the space of
matrices and some information easier to get in this setting. It then
could be (partially) lifted because  of 1-1 correspondence of
orbits. More precisely, it could be lifted in sense of
parametrization, but if we consider, for example,  degeneration of
orbits situation may changes after restriction of them.

In this section we will give a parametrization (by two parameters)
of the family ${\cal M}_n$ of representations defined by ${\rm rk} Y
= n-1$.
What we actually doing here, we obtain this parametrization for
$W_{(n)}$. Due to 1-1 correspondence between the orbits we then have
a parametrization of ${\cal M}_n$.

Let we restrict the orbits even a little further, considering the
action of the group $G = Z(J_{\cal P}) \cap SL_n$, where the 1-1
correspondence with the initial orbits will be clearly preserved. In
the case ${\cal P}= (n)$ the group $G$ can be presented as follows:
$$G=\{I+\a_1 Y + \a_2 Y^2 + ... + \a_{n-1} Y^{n-1}\},$$
\noindent due to our description of the centralizer of $Y$ in
section \ref{descript}. This group acts on the affine space of the
dimension $n$:
$$W_{(n)}=\{\l I + X^0+c_1 Y + c_2 Y^2 + ... + c_{n-1} Y^{n-1}\}$$
\noindent here $\l$ is the eigenvalue of $X$ and $X^0$ is the matrix
defined in section \ref{descript} with the second diagonal $[0,1,
\dots,(n-1)]$ and zeros elsewhere.


Let we fix first the eigenvalue: $\, \l=0$, we get then the space of
dimension $n-1$:
$$W_{(n)}'=\{ X^0+c_1 Y + c_2 Y^2 + ... + c_{n-1}
Y^{n-1}\}.$$

We intend to calculate now the dimension of the orbit ${\cal
O}_{(n)}(X,G)$ of $X$ with fixed eigenvalue $\l=0$ under $G$ --
action.

Consider the map $\phi: G \longrightarrow W_{Y}'$ defined by this
action: $\phi (C) = CXC^{-1}$, then ${\rm Im}\, \phi ={ \cal
O}_{(n)}(X,G)$. We are going to  calculate the rank of Jacobian of
this map. We will see that it is constant on $G$ and equals to
$n-2$. This tells us that each orbit ${ \cal O}_{(n)}(X,G)$ is an
$n-2$ dimensional manifold and hence there couldn't be more then 2
parameters involved in parametrization of orbits.

\subsection{Calculation of the rank of Jacobian}

\begin{theorem}\label{t7.1} Let $G$ be an intersection of $SL_n$
with the centralizer of $Y$. Consider the action of this group on
the affine space $W_{Y}'=\{ X^0+c_1 Y + c_2 Y^2 + ... + c_{n-1}
Y^{n-1}\}$ by conjugation. Then the rank of the Jacobian of the map
$\phi: G \longrightarrow W_{Y}'$ is equal to $n-2$ in any point $C
\in G$.
\end{theorem}

\proof
Consider $d \phi(C)(\D)=(C+\D)^{-1}X(C+\D)-C^{-1}XC,$ where

$$C=I+\a_1 Y + \a_2 Y^2 + ... + \a_{n-1}
Y^{n-1},$$

$$X=X^0+c_1 Y + c_2 Y^2 + ... + c_{n-1} Y^{n-1},$$

$$\D=\b_1 Y + \b_2 Y^2 + ... + \b_{n-1} Y^{n-1}.$$

Let we present $(C+\D)^{-1}$ in the following way:

$$(C+\D)^{-1}=(I+\D C^{-1})^{-1} C^{-1}=$$

$$(I-\D C^{-1} +  {\rm lower \, order \, terms \, on}
\D) C^{-1}.$$

Then

$$(C+\D)^{-1}X(C+\D)-C^{-1} X C =  $$

$$(I-\D C^{-1} +
{\rm lower \, order \, terms \, on}  \D) C^{-1} X (C+\D) - C^{-1} X
C=$$

$$ - \D C^{-2} X C + C^{-1} X \D + {\rm lower \, order \, terms \, on}\D =$$

$$(-\D C^{-1} \cdot C^{-1} X + C^{-1} X \cdot \D C^{-1}) C
+ {\rm lower \, order \,terms \, on}\D .$$

Denote  $\tilde \D := \D C^{-1}$ and $\tilde X :=  C^{-1} X$.
Obviously multiplication by $C$ preserves the rank and rank of
linear map $d \phi(C)(\D)$ is equal to the rank of the map $T(\tilde
\D) = [\tilde X, \tilde \D]$.

Here again
$\tilde \D$ has
a form
$$\tilde \D=
\g_1 Y + \g_2 Y^2 + ... + \g_{n-1} Y^{n_1}.$$

Let us compute commutator of $\tilde X$ with $Y^k$. Taking into
account that $C^{-1}$ is a polynomial on $Y$, hence commute with
$Y^k$ and also the relation in  algebra $R$: $ XY^k-Y^kX=k Y^{k+1}.$
We get $ \tilde X Y - Y \tilde X=
C^{-1}XY^k-Y^kC^{-1}X=C^{-1}(XY^k-Y^kX)= C^{-1} k Y^{k+1} $. Hence
$$ \tilde X p(Y) - p(Y) \tilde X = C^{-1} Y^2 p'(Y)$$
\noindent for arbitrary polynomial $p$. Applying this for the
polynomial $\tilde \D$ we get
$$ T(\tilde \D)= [\tilde X, \tilde \D] =
\sum_{k=1}^{n-2} \g_k k C^{-1} Y^{k+1},$$ \noindent hence this
linear map has rank $n-2$. $\Box$






From the theorem \ref{t7.1} we could deduce the statement concerning
parametrization of isoclasses of modules in the family ${\cal
M}_{n}$ .

We mean by {\it parametrization} (by $m$ parameters) the existence
of $m$ smooth algebraically independent functions which are constant
on the orbits and separate them.

\begin{corollary}\label{c7.2}
Let ${\cal U}_{(n)}$ be the stratum as above. Then the set of
isomorphism classes of indecomposable modules from ${\cal U}_{(n)}$
could be parameterized by at most two parameters.
\end{corollary}

\proof Directly from the theorem \ref{t7.1} applying the theorem on
locally flat map \cite{DNF} to $\phi: G \longrightarrow W_{(n)}'$ we
have that ${\rm Im} \phi ={\cal O}_{(n)}(X,G)$ is an $n-2$
dimensional manifold. We have to mention here that this is due to
the fact that the image has no selfintersections. This is the case
since the preimage of any point $P$ is connected (it is formed just
by the solutions of the equation $CX=PC$ for $C \in G$). Hence we
can parametrize these orbits lying in the space $W_{(n)}$ of
dimension $n$ by at most two parameters. Due to 1-1 correspondence
to the whole orbits ${\cal O}(X,Y)$ the latter have the same
property.$\Box$

\begin{proposition}\label{p7.1} Parameters $\mu$ and $\lambda$ are invariant under the
action of $G$ on the set of matrices $\displaystyle \tiny
\left\{\left(\begin{array}{cccccc}
\lambda&\mu+1&&&&\\
&\lambda&\mu+2&&\smash{\hbox{\normalsize*}}&\\
&&\lambda&\mu+3&&\\
&&&\ddots&\ddots&\\
&\smash{\hbox{\normalsize0}}&&&\lambda&\mu+n-1\\
&&&&&\lambda\end{array}\right)\right\}.$
\end{proposition}

\proof Direct calculation of $ZMZ^{-1}$ for $Z\in G$ as described
above shows that elements in first two diagonals of $M$ will be
preserved. $\Box$

Hence from the corollary \ref{c7.2} and proposition \ref{p7.1} we
have the following classification result for the family ${\cal M}_n$
of representations with full Jordan block $Y$, or equivalently with
the condition $n - rk Y = 1$.

\begin{theorem}\label{t7.2} Let $P_{\lambda, \mu}$ denotes the pair
$(X_{\lambda, \mu},Y)$, where
$$\tiny X_{\lambda,\mu}=\left(\begin{array}{cccccc}
\lambda&\mu+1&&&&\\
&\lambda&\mu+2&&\smash{\hbox{\normalsize0}}&\\
&&\lambda&\mu+3&&\\
&&&\ddots&\ddots&\\
&\smash{\hbox{\normalsize0}}&&&\lambda&\mu+n-1\\
&&&&&\lambda\end{array}\right), \ \ \
Y=\left(\begin{array}{cccccc}
0&1&&&&\\
&0&1&&\smash{\hbox{\normalsize0}}&\\
&&0&1&&\\
&&&\ddots&\ddots&\\
&\smash{\hbox{\normalsize0}}&&&0&1\\
&&&&&0\end{array}\right).
$$
Every pair $(X,Y) \in {\cal M}_n$ is conjugate to $P_{\lambda, \mu}$
for some $\lambda, \mu$. No two pairs $P_{\lambda, \mu}$ with
different $(\lambda, \mu) $ are conjugate.

\end{theorem}

Let we mention that number of parameters does not depends of $n$ in
this case.





\subsection{Some examples of tame strata (up to auto-equivalence)}

We collect (quite rare) examples of tame strata in the suggested
above stratification related to the Jordan normal form of $Y$. We
present here tameness results for the representation type of
families of reps lying in the stratum ${\cal U}_{(n-1,1)}$ with
respect to auto-equivalence relation on modules. It was defined in
section \ref{rfd} and consists of gluing orbits which could be
obtained one from another using automorphism of the initial algebra.

\begin{theorem}\label{t7.4}
The subset  of all n-dimensional representations corresponding to
$Y$ with full Jordan block, or equivalently defined by the condition
$n - rk Y=1$, has a finite representation type with respect to
auto-equivalence relation on modules.

The subset  of all n-dimensional representations corresponding to
$Y$ with the Jordan structure ${\cal P}=(n-1,1)$ is tame, that is
parametrizable by one parameter, with respect to auto-equivalence
relation. \end{theorem}

\proof The proof analogues to the proof of the theorem \ref{t7.2}.
We present here pictures showing how $X$ and $Y$ act on basis and
where parameters appear:
$$
\begin{array}{ll}
{\cal P}=(n)&\\
&\bull{e_1}\doubarr{y}{x}\bull{e_2}\doubarr{y}{(x,2)}\bull{e_3} \
\ {\dots} \  \ \bull{e_{n-1}}\doubarr{y}{(x,n-1)}\bull{e_n}\\
\phantom0&\\
\phantom0&\\
\phantom0&\\
{\cal P}=(n-1,1)&\\
&\bull{e_1}\doubarr{y}{x}\bull{e_2}\doubarr{y}{(x,2)}\bull{e_3} \ \
{\dots} \  \ \bull{e_{n-2}}\doubarr{y}{(x,n-2)}\bull{e_{n-1}}
\siglearr{(x,\alpha)}\bull{e_{n}}\\ &\uparrow\hskip7.25cm|\\
&\,\hskip.47mm\smash{\raise 9.5pt\hbox to 7.47cm{\hrulefill}}\\
&\hfill\smash{\raise 10pt\hbox{$\scriptstyle
(x,\alpha^{-1})$}}\hfill
\end{array}
$$$\Box$

\section{The case of one block and Ringel's classification of
complete local algebras}\label{ring}

As we have shown above the set of orbits corresponding to the
full-block Jordan structure of $Y$ in the variety of n-dimensional
modules could be parametrized by two parameters. Therefore this
family  ${\cal M}_n $ of representations defined in section
\ref{quiv} is wild. Nevertheless we intend to prove here that all
representations from  ${\cal M}_n $ have only one finite dimensional
algebra (for any dimension $n$) as their image.

We shall show the place of these algebras in Ringel's classification
of complete local algebras \cite{Ringloc} by calculating their
defining relations. Let we remind that as we have seen in the
Corollary \ref{cindlocal} any indecomposable representation has a
local algebra as an image, particularly, representations with
full-block $Y$ do. We are going to prove here that for $n \leq 4$
all image algebras $A_n$ are tame and for $n \geq 5$ they are wild.

\begin{theorem}\label{tonei} Let $\rho_n:R\to M_n(k)$ be a finite dimensional
representation of $R$, where $Y=\rho(y)$ has a full-block Jordan
structure. Then the image algebra $A_n=\rho_n(R)$ does not depend on
the choice of $\rho_n$.
\end{theorem}

\proof In order to calculate the linear basis of $A_n$ it suffices
to find matrices $\epsilon_n(y^nx^m)$, which are images of normal
monomials $y^nx^m$ under the representation $\epsilon_n$ defined in
section \ref{rfd}. It is enough to consider $\epsilon_n$ since as we
have shown in section  \ref{rfd}, any $\rho_n$ is auto-equivalent to
$\epsilon_n$ and the images of auto-equivalent representations
coincide. These matrices were calculated in the lemma \ref{lImyx}.
Namely, the matrix $\epsilon_n(y^{k}x^m)$ has the vector
$p(o),p(1),...,p(j),...$ of values of  polynomial
$p(j)=(k+j)...(k+m+j-1)=\prod_{i=1}^m (k+j+i)$ in the  diagonal
number $k+m$ and zeros elsewhere. The linear span of these matrices
gives us the desired image algebra. $\Box$

Recall that a $k$-algebra $A$ is called {\it local} if $A=k\oplus
{\rm Jac}(A)$, where ${\rm Jac}(A)$ is the Jacobson radical of $A$.
One can also consider the {\it completion} of $A$:
$\overline{A}=\lim\limits_{\displaystyle\longleftarrow}A/({\rm
Jac}(A))^n$. An algebra $A$ is called complete if $A=\overline A$.

It was shown in section \ref{descript} that in the case of
full-block $Y$, $X$ has only one eigenvalue. We also have proved in
Theorem \ref{t2.3} that for all image algebras  their semisimple
part $A/{\rm Jac}(A)$ is the direct sum of $r$ copies of $k$, $r$
being the number of different eigenvalues of $X$. Hence in the
full-block case the image algebra is local. It is also complete,
because ${\rm Jac}(A)^N=0$ for $N$ large enough. Indeed, we can use
here the Corollary \ref{radft} which describe the radical, or
observe directly that since $A=k\oplus {\rm Jac}(A)$ and $A$
consists of polynomials on $X^0$ and $J_n$ (as an image of one of
representations $\epsilon: (x,y) \mapsto (X^0,J_n)$), then $Jac(A)$
consists of those polynomials which have no constant term. Since the
matrices $J_n$ and $X^0$ are nilpotent of degree $n$ and $n-1$
respectively, ${\rm Jac}(A)^{2n}=0$.

\begin{theorem}  The image algebra $A_n$ of a representation
$\rho_n \in {\cal M}_n $ is wild for any $n\geq 5$. It has a
quotient isomorphic to the wild algebra given by relations
$y^2,yx-xy,x^2y,x^3$ from the Ringel's list of minimal wild local
complete algebras. The image algebras $A_1,A_2$ and $A_3$ are tame.
\end{theorem}

\proof We intend to show that for $n$ big enough, the algebra $A_n$
has a quotient isomorphic to the algebra $W=\langle
x,y|y^2=yx-xy=x^2y=x^3=0 \rangle$, which is number c) in the
Ringel's list of minimal wild local complete algebras
\cite{Ringloc}.

The algebra $W$ is 5-dimensional. Let us consider the ideal $J$ in
$A_n$ generated by the relations above on the image matrices $X$ and
$Y$. This ideal has obviously codimension not exceeding 5. We intend
to show that $J$ has codimension exactly 5 and therefore $W$ should
be isomorphic to $A/J$.

Let us look at the ideal $J$, which is generated by
$\{Y^2,X^2Y,X^3,XY-YX\}$. First, since $XY-YX=Y^2$, $J$ is generated
by $\{Y^2,X^2Y,X^3\}$. It is easy to see that $Y^2$ has zeros on
first two diagonals and the vector ${\bf {1}}=(1,\dots,1)$ on the
third one,  $X^2Y$, $X^3$ have zeros on the first three diagonals.
An arbitrary  element of $A_n$ has the constant vector $c{\bf 1}$
for some $c \in k$ on the main diagonal.
Hence we see that a general element of the ideal $J$ has zeros on
the first two diagonals and the constant sequence $c{\bf 1}$ on the
third one. Taking into account that $A_n$ comprises the upper
triangular matrices that have values of a polynomial of degree at
most $m$ on $m$-th diagonal (lemma \ref{lImyx}), we see that the
main diagonal gives an impact of $1$ to the codimension of $J$, the
first diagonal gives an impact of $2$ to the codimension of $J$ if
the length of this diagonal is at least $2$ (that is $n\geq 3$) and
the second diagonal, ---  an impact of $2$, if the length of this
diagonal is at least $3$ (that is $n\geq 5$). Thus, dim$\,A/J\geq 5$
if $n\geq 5$. This completes the proof in the case $n\geq 5$.

Tameness of $A_1$ and $A_2$ is obvious. For $A_3$ the statement
follows from the dimension reason: dim$\,A_3=4$, it is less then the
dimensions of all 2-generated algebras from the Ringel's list of
minimal  wild algebras. Since his theorem (theorem 1.4 in
\cite{Ringloc}) states that any local complete algebra is either
tame of has a quotient from the list, $A_3$ can not be wild by
dimension reasons. Hence $A_3$ is tame. $\Box$

Let us consider now the case $n=4$.

\begin{theorem}  Let $\rho_4 \in {\cal M}_4$ be a four dimensional
representation of the algebra $R$. Then the image algebra
$A_4=\rho_4 (R)$ is given by the relations $k\langle x,y|x^2=-2xy,
xy=yx+y^2,x^3=0\rangle$ and is tame.
\end{theorem}

\proof We intend to show that no one of the algebras from the
Ringel's list of minimal wild algebras can be obtained as a quotient
of $A_4$. After that using the Ringel's theorem, we will be able to
conclude that it is tame. Suppose that there exists an ideal $I$ of
$A_4$ such that $A_4/I$ is isomorphic to $W_j$ for some $j=1,2,3,4$,
where
\begin{align*}
W_1&=k\langle u,v|u^2,uv-\mu vu\ (\mu\neq0),v^2u,v^3\rangle,
\\
W_2&=k\langle u,v|u^2,uv,v^2u,v^3\rangle,
\\
W_3&=k\langle u,v|u^2,vu,uv^2,v^3\rangle,
\\
W_4&=k\langle u,v|u^2-v^2,vu\rangle.
\end{align*}

Since all $W_j$ are 5-dimensional and $A_4$ is 6-dimensional, the
ideal $I$ should be one-dimensional. Due to our knowledge on the
matrix structure of the algebra $A_4$, we can see that there is only
one one-dimensional ideal $I_4$ in $A_4$ and that $I_4$ consists of
the matrices with at most one non-zero entry being in the upper
right corner of the matrix:
$$
I_4=\left\{\left(
\begin{array}{cccc}
0&0&0&c\\0&0&0&0\\0&0&0&0\\0&0&0&0
\end{array}
\right)\right\}.
$$
After factorization by this ideal we get a $5$-dimensional algebra
given by relations
$$
\overline{A}_4=A_4/I=k\langle x,y|x^2=-2xy,
xy=yx+y^2,x^3=0,y^3=0\rangle.
$$
The question now is whether this algebra is isomorphic to one of the
algebras from the above list. Suppose that there exists an
isomorphism $\phi_j:W_j\to \overline{A}_4$ for some
$j\in\{1,2,3,4\}$. Denote $\phi_j(u)=f_j$ and $\phi_j(v)=g_j$.

First, let us mention that $f_j$ and $g_j$ have zero free terms:
$f_j^{(0)}=g_j^{(0)}=0$ because the equalities $\phi_j(u^2)=f_j^2=0$
and $\phi_j(v^3)=g_j^3$ imply $(f_j^{(0)})^2=(g_j^{(0)})^3=0$ and
therefore $f_j^{(0)}=g_j^{(0)}=0$ if $j=1,2,3$ and the equalities
$\phi_4(u^2-v^2)=f_4^2-g_4^2=0$ and $\phi_4(uv)=f_4g_4=0$ imply
$(f_4^{(0)})^2=(g_4^{(0)})^2$ and $f_4^{(0)}g_4^{(0)}=0$ and
therefore $f_4^{(0)}=g_4^{(0)}=0$.

The second observation is that the terms of degree 3 and more are
zero in $A_4$. Therefore we can present the polynomials $f_j$ and
$g_j$ as the sum of their linear and quadratic (on $x$ and $y$)
parts. So, let $f_j=f_j^{(1)}+f_j^{(2)}$ and
$g_j=g_j^{(1)}+g_j^{(2)}$, where
$$
f_j^{(1)}=ax+by,\ \ g_j^{(1)}=\alpha x+\beta y,\ \
f_j^{(2)}=cyx+dy^2,\ \ g_j^{(2)}=\gamma yx+\delta y^2.
$$
In order to get entire  linear part of the algebra $A_4$ in the
range of $\phi_j$ we need to have
\begin{equation}
{\rm det}\,\left|
\begin{array}{cc}
a&b\\\alpha&\beta
\end{array}
\right|\neq0. \label{det}
\end{equation}

For any $j=1,2,3,4$ we are going to obtain a contradiction of the
last condition with the equations on $a,b,\alpha,\beta$ coming from
the relations of the algebra $W_j$.

For instance, consider the case $j=2$. From
$0=u^2=f_j^2=(f_j^{(1)})^2=2(ab-a^2)yx+(b^2+ab-2a^2)y^2$ we get
$2a(b-a)=0$ and $b^2+ab-2a^2=0$. The first equation gives us that
either $a=0$ or $a=b$. In the case $a=0$ the second equation implies
$b=0$ and the equality $a=b=0$ already contradicts (\ref{det}).
Another solution is $a=b\neq 0$. From
$0=uv=f_jg_j=f^{(1)}g^{(1)}=(ax+by)(\alpha x+\beta y)$, substituting
$a=b$, we get $0=a(x+y)(\alpha x+\beta y)=a(\beta-\alpha)(yx+2y^2)$.
Hence $\beta=\alpha$, which together with the equality $a=b$
contradicts (\ref{det}).

In the other three cases one can  get a contradiction with
(\ref{det}) along the same lines, which completes the proof. $\Box$




\section{Acknowledgments}

This work was carried out in 2003-2004 during my stay at the
Max-Planck-Intit\"ut f\"ur Mathematik in Bonn. I am thankful to this
institution for hospitality and excellent research atmosphere and to
many collegues with whom I have been discussing ideas related to
this paper.

This work also enjoys a partial support of the Pythagoras II
programm of Eurocomission.

\end{document}